\documentclass[1p]{amsart} 
\usepackage{mathptmx,amsmath,mathbbol,amssymb}

\newtheorem{theorem}{Theorem}[section]
\newtheorem{proposition}[theorem]{Proposition}
\newtheorem{lemma}[theorem]{Lemma}
\newtheorem{corollary}[theorem]{Corollary}
\newtheorem{example}[theorem]{Example}
\newtheorem{remark}[theorem]{Remark}
\newtheorem{definition}[theorem]{Definition}

\newcommand{\klg}[1]{\left\{#1\right\}}
\newcommand{\klr}[1]{\left(#1\right)}
\newcommand{\kle}[1]{\left[#1\right]}
\newcommand{\klrk}[1]{(#1)}
\newcommand{\klek}[1]{[#1]}
\newcommand{\klgk}[1]{\{#1\}}
\newcommand{\betrag}[1]{{\left|#1\right|}}
\newcommand{\betragk}[1]{{|#1|}}
\newcommand{\norm}[1]{\left\|#1\right\|}
\newcommand{\normk}[1]{\|#1\|}
\newcommand{\normka}[1]{\|\;#1\;\|}

\def\nt{\triangledown_2}
\def\dt{\triangle_2}

\def\supp{{\rm supp}}

\def\loc{{\rm loc}}
\def\IR{{\mathbb{R}}}
\def\IN{{\mathbb{N}}}

\def\Cap{\operatorname{Cap}}

\def\domain{\mathsf{D}}

\def\cB{\mathcal{B}}
\def\cD{\mathcal{D}}
\def\cE{\mathcal{E}}
\def\cX{\mathcal{X}}
\def\cY{\mathcal{Y}}

\def\cK{\mathcal{K}}

\def\cN{\mathcal{N}}
\def\cP{\mathcal{P}}

\def\sC{\;\mathsf{C}}

\def\sV{\mathcal{V}}

\def\tW{\tilde W}

\def\se{\mathsf{e}}
\def\su{\mathsf{u}}

\def\sv{\mathsf{v}}
\def\sw{\mathsf{w}}

\def\prmz{(M1)}
\def\prme{(M2)}
\def\prmx{(M3)}
\def\prmc{(L1)'}

\def\prtl{(V1)}
\def\prd{(V2)}
\def\prss{(V3)}
\def\prvl{(V4)}
\def\prn{(V5)}
\def\prnz{(V6)}

\def\prchm{(C1)}
\def\prchi{(C2)}
\def\prchd{(C3)}
\def\prtcc{(L1)}
\def\prtnr{(L2)}
\def\prtrn{(L3)}
\def\prco{(K1)}
\def\prcc{(K2)}
\def\prcr{(K3)}
\def\prxm{(X1)}
\def\prxe{(X2)}
\def\prmh{(M1)}
\def\prmo{(M2)}
\def\prmz{(M3)}
\def\prmc{(M4)}

\begin{document}

\title{On a Capacity for Modular Spaces}
\author{Markus Biegert}
\email{markus.biegert@uni-ulm.de}
\address[rvt]{Institute of Applied Analysis, University of Ulm, 89069 Ulm, Germany}

\begin{abstract}
  The purpose of this article is to define a capacity on certain topological measure 
  spaces $X$ with respect to certain function spaces $V$ consisting of measurable functions.
  In this general theory we will not fix the space $V$ but we emphasize that
  $V$ can be the classical Sobolev space $W^{1,p}(\Omega)$,
  the classical Orlicz-Sobolev space $W^{1,\Phi}(\Omega)$,
  the Haj{\l}asz-Sobolev space $M^{1,p}(\Omega)$, the Musielak-Orlicz-Sobolev space
  (or generalized Orlicz-Sobolev space) and many other spaces. Of particular
  interest is the space $V:=\tW^{1,p}(\Omega)$ given as the closure of
  $W^{1,p}(\Omega)\cap C_c(\overline\Omega)$ in $W^{1,p}(\Omega)$.
  In this case every function $u\in V$ (a priori defined only on $\Omega$) has a trace
  on the boundary $\partial\Omega$ which is unique up to a $\Cap_{p,\Omega}$-polar set.
\end{abstract}

\maketitle


\section{Introduction}

The notion of capacity is fundamental to the analysis of pointwise behavior of Sobolev type
functions. Depending on the starting point of the study, the capacity of a set can be
defined in many appropriate ways. The Choquet theory \cite{choquet:54:toc} gives a
standard approach to capacities. Capacity is a necessary tool in classical and
nonlinear potential theory. One purpose of this article is to introduce an extension of the
classical $p$-capacity which we call the {\em relative $p$-capacity}.
For example, given an open set $\Omega\subset\IR^N$ the
classical $p$-capacity and the relative $p$-capacity can be used to decide whether a
given function $\su$ lies in $W^{1,p}_0(\Omega)$ or not.
The notion of relative $2$-capacity was
first introduced by Wolfgang Arendt and Mahmadi Warma in \cite{arendt:03:lrb} to study
the Laplacian with general Robin boundary conditions on arbitrary domains in $\IR^N$.
For the investigation of the $p$-Laplacian with generalized Robin boundary conditions on bad
domains, such as the snowflake (the domain bounded by the von Koch curve), the relative
$p$-capacity plays an important role.

For results on the classical $p$-capacity and other capacities we refer the
reader to the following books and the references therein:
David R. Adams and Lars I. Hedberg \cite{adams:96:fsp},
Nicolas Bouleau and Francis Hirsch \cite{bouleau:91:df},
Gustave Choquet \cite{choquet:54:toc},
Lawrence C. Evans and Ronald F. Gariepy \cite{evans:92:mtf},
Juha Heinonen and Tero Kilpel{\"a}inen and Olli Martio \cite{heinonen:93:npt},
Jan Mal\'y and William P. Ziemer \cite{ziemer:97:fr} and
Vladimir G. Maz'ya \cite{mazya:85:ssp}.
For capacities on Orlicz-Sobolev spaces and their fine behavior we refer to a recent
article of J.~Mal\'y, D.~Swanson and W.~P.~Ziemer \cite{ma:09:fbfos} and the references therein.


\section{Preliminaries}

\subsection{Setting}

\begin{definition}[Type $\Lambda$]
  Let $V$ be a real vector space. We will call a mapping $\rho:V\to [0,\infty]$ a
  {\em Luxemburg functional} on $V$ and $V=(V,\rho)$ a {\em Luxemburg space} if
  \begin{itemize}
    \item[\prmz] $\rho(\su)=0$ if and only if $\su=0$;
    \item[\prme] $\rho(-\su)=\rho(\su)$ for all $\su\in V$;
    \item[\prmx] $\rho(\alpha\su+\beta\sv)\leq \alpha\rho(\su)+\beta\rho(\sv)$ whenever
                 $\alpha+\beta=1$, $\alpha,\beta\geq 0$ and $\su,\sv\in V$;
    \item[\prmc] $\lim_{\lambda\to 0+} \rho(\lambda\su)=0$ for all $\su\in V$;
  \end{itemize}
  A Luxemburg space $(V,\rho)$ is endowed with the {\em Luxemburg norm} $\norm{\cdot}_V$ 
  given by
  \[ \norm{\su}_{V} := \inf\klg{\alpha>0:\rho(\su/\alpha)\leq 1}.
  \]
  See \cite[Theorem 1.5]{mu:83:osms}. If in addition
  \begin{itemize}
    \item[\prtcc] $\rho:(V,\norm{\cdot}_V)\to\IR_+:=[0,\infty)$ is continuous;
    \item[\prtnr] $\lim_{\norm{\su}_V\to\infty} \rho(\su)=\infty$;
    \item[\prtrn] $\lim_{\rho(\su)\to 0} \norm{\su}_V=0$;
  \end{itemize}
  then we say that $(V,\rho)$ is of type $\Lambda$. A Luxemburg space $(V,\rho)$ is called
  reflexive/complete if $V$ with respect to the Luxemburg norm is reflexive/complete.
\end{definition}

\begin{remark}\label{rem:semicont}
  If $(V,\rho)$ is a Luxemburg space, then $\rho(\su)\leq\norm{\su}_V$ for all
  $\su\in V$ with $\norm{\su}_V<1$ \cite[Theorem 1.5(III)]{mu:83:osms}.
  If in addition $\rho$ is continuous [see property \prtcc], then $\rho(\su/\norm{\su}_V)=1$
  for all $\su\in V\setminus\klg{0}$. Moreover, $\rho$ is weakly lower-semicontinuous
  \cite[Corollary III.8]{brezis:83:af}, in particular,
  \[ \rho(\su) \leq \liminf_{\su_n\rightharpoonup\su} \rho(\su_n).
  \]
  If $(V,\norm{\cdot}_1)$ is a real normed vector space, then $\rho(\su):=\norm{\su}_1$
  is a Luxemburg functional on $V$, $(V,\rho)$ is a Luxemburg space of type $\Lambda$
  and the Luxemburg norm $\norm{\cdot}_V$ is equal to $\norm{\cdot}_1$. So completeness and
  reflexivity of a Luxemburg space of type $\Lambda$ is not automatic.
\end{remark}

\begin{definition}[Type $\Theta$]
  By a {\em topological measure space} (abbreviated by {\em tms}) we mean the quadruple $(X,\tau,\Sigma,\mu)$
  where $(X,\tau)$ is a topological space and $(X,\Sigma,\mu)$ a measure space.
  We will say that a topological measure space $X$ is of {\em type $\Theta$} if
  \begin{itemize}
    \item[\prmh] $(X,\tau)$ is a Hausdorff topological space;
    \item[\prmo] $(X,\Sigma,\mu)$ is a complete measure space;
    \item[\prmz] If $U\subset X$ is open and $\Sigma$-measurable with $\mu(U)=0$, then $U=\emptyset$;
    \item[\prmc] Every open set $U\subset X$ is the countable union of compact sets.
  \end{itemize}
\end{definition}

\begin{remark}
  If $X=(X,\tau,\Sigma,\mu)$ is of type $\Theta$, $f,g\in C(X)$ are measurable and
  $f=g$ $\mu$-a.e.\ on $X$, then $f=g$ everywhere on $X$.
\end{remark}

\begin{definition}[Domination of type $\Theta$]
  For $j=1,2$ let $X_j=(X_j,\tau_j,\Sigma_j,\mu_j)$ be tms of type $\Theta$.
  We will say that $X_2$ dominates $X_1$, abbreviated by $X_1\preceq X_2$, if
  \begin{itemize}
    \item[(D1)] $X_1\in\Sigma_2$ and $\Sigma_1=\Sigma_2\cap X_1=\klg{M\cap X_1:M\in\Sigma_2}$;
    \item[(D2)] $\tau_1=\tau_2\cap X_1=\klg{O_2\cap X_1:O_2\in\tau_2}$;
    \item[(D3)] $\mu_1(A)\leq\mu_2(A)$ for all $A\in\Sigma_1$.
  \end{itemize}
\end{definition}

\begin{example}\label{ex:tms}
  Let $\Omega_1\subset\Omega_2$ be two non-empty open sets in $\IR^N$,
  $\lambda$ be the $N$-dimensional Lebesgue measure and $\Sigma$ be the
  $\sigma$-algebra of all Lebesgue-measurable sets in $\IR^N$.
  \begin{enumerate}
    \item For $j=1,2$ we let $X_j:=\Omega_j$, $\tau_j:=\tau_{\IR^N}\cap X_j$,
          $\Sigma_j:=\Sigma\cap X_j$ and 
          define the measure $\mu_j$ on $(X_j,\Sigma_j)$ by $\mu_j(B):=\lambda(B\cap X_j)$.
          Then $(X_j,\tau_j,\Sigma_j,\mu_j)$ is a tms of type $\Theta$ and $X_1\preceq X_2$.
    \item For $j=1,2$ we let $X_j:=\overline\Omega_j$, $\tau_j:=\tau_{\IR^N}\cap X_j$,
          $\Sigma_j:=\Sigma\cap X_j$ and 
          define the measure $\mu_j$ on $(X_j,\Sigma_j)$ by $\mu_j(B):=\lambda(B\cap X_j)$.
          Then $(X_j,\tau_j,\Sigma_j,\mu_j)$ is a tms of type $\Theta$ and $X_1\preceq X_2$.
  \end{enumerate}
\end{example}

\begin{definition}[Class $\Upsilon$]
  Let $X=(X,\tau,\Sigma,\mu)$ be a tms of type $\Theta$ and denote by $L^0(X)$ the vector space of all
  real-valued (equivalence classes of $\mu$-a.e.\ equal) measurable functions on $X$.
  A subspace $V\subset L^0(X)$ equipped with a Luxemburg functional $\rho$ belongs
  to the {\em class $\Upsilon=\Upsilon(X,\tau,\Sigma,\mu)$}, briefly $V\in\Upsilon$,
  if it satisfies the following properties.
  \begin{itemize}
    \item[\prtl]  $(V,\rho)$ is a reflexive and complete Luxemburg space of type $\Lambda$;
    \item[\prd]   The space $V\cap C_c(X)$ is dense in $V$;
    \item[\prss]  If $\su_n,\su\in V$ and $\su_n\to\su$ then
                  a subsequence of $(\su_n)_n$ converges $\mu$-a.e.\ to $\su$;
    \item[\prvl]  $V$ is a vector lattice with respect to the $\mu$-a.e.\ pointwise ordering;
    \item[\prn]   The function $\su\wedge c$ belongs to $V$ for every $\su\in V$ and $c\in\IR_+$
                  and $\rho(\su\wedge c)\leq\rho(\su)$;
    \item[\prnz]  For every $c\in\IR_+$ the mapping $V\to V$, $\su\mapsto\su\wedge c$ is continuous.
  \end{itemize}
\end{definition}

\begin{remark}
  From property \prn\ we get that $\rho(\su^+)\leq\rho(\su)$ for all $\su\in V$.
  This implies that $\norm{\su^+}_V\leq\norm{\su}_V$. Similarly, we get that
  $\rho(\su^-)=\rho((-\su)^+)\leq \rho(-\su)=\rho(\su)$ and $\normk{\su^-}_V\leq\normk{\su}_V$.
  Therefore
  \[ \normka{\betrag{\su}}_V=\normk{\su^++\su^-}_V\leq
     \normk{\su^+}_V+\normk{\su^-}_V \leq \normk{\su}_V+\normk{\su}_V =
     2\normk{\su}_V.
  \]
\end{remark}

\begin{definition}[Domination of class $\Upsilon$]\label{def:domups}
  For $j=1,2$ let $X_j$ be a tms of type $\Theta$ and
  $(V_j,\rho_j)$ of class $\Upsilon(X_j)$. Then we say that {\em $V_2$ dominates $V_1$},
  abbreviated by $V_1\preceq V_2$, if
  $X_1\preceq X_2$,
  $\su_2|_{X_1}\in V_1$ for all $\su_2\in V_2$ and there is a constant $c>0$ such that
  $\rho_1(\su_2|_{X_1})\leq c\rho_2(\su_2)$ for all $\su_2\in V_2$.
\end{definition}

\begin{example}[Sobolev spaces]\label{ex:sobo-one}
  For $p\in(1,\infty)$ and $\Omega\subset\IR^N$ open we let $W^{1,p}(\Omega)\subset L^p(\Omega)$
  be the first order {\em Sobolev space} consisting of all functions $\su\in L^p(\Omega)$
  whose distributional derivatives of order one belong to $L^p(\Omega)$. Equipped with
  the norm $\norm{\cdot}_{W^{1,p}(\Omega)}$ given by
  \[ \norm{\su}^p_{W^{1,p}(\Omega)}:=
     \normk{\su}^p_{L^p(\Omega)} + \normka{\betrag{\nabla\su}}^p_{L^p(\Omega)}
  \]
  the space $W^{1,p}(\Omega)$ is a reflexive Banach space. Let $V:=\tW^{1,p}(\Omega)$ be the
  closure of $W^{1,p}(\Omega)\cap C_c(\overline\Omega)$ in $W^{1,p}(\Omega)$ and
  $\rho_{1,p,\Omega}(\su):=\norm{\su}_{W^{1,p}(\Omega)}$.
  Then $(V,\rho_{1,p,\Omega})$ is of class $\Upsilon(\overline\Omega)$ where the tms
  $X:=\overline\Omega$ is as in Example \ref{ex:tms}(2).
  Moreover, if $\Omega_1\subset\Omega_2$ are non-empty open sets in $\IR^N$
  we get that 
  $\tW^{1,p}(\Omega_1)\preceq\tW^{1,p}(\Omega_2)$.
\end{example}

\begin{definition}[$\cN$-function]
  A mapping $\Phi:\IR\to\IR_+$ is called an $\cN$-function if
  \begin{itemize}
    \item[(N1)] $\Phi$ is even and convex;
    \item[(N2)] $\Phi(x)=0$ if and only if $x=0$;
    \item[(N3)] $\lim\limits_{x\to 0+} x^{-1}\Phi(x)=0$ and 
                $\lim\limits_{x\to\infty} x^{-1}\Phi(x)=\infty$.
  \end{itemize}
  Let $\Psi:\IR\to\IR_+$ be given by $\Psi(y):=\sup\klg{x\betrag{y}-\Phi(x):x\geq 0}$.
  Then $\Psi$ is an $\cN$-function, called the {\em complementary $\cN$-function} to $\Phi$.
\end{definition}

\begin{definition}[The $\dt$- and $\nt$-condition]
  An $\cN$-function $\Phi$ is said to obey the global {\em $\dt$-condition} if
  there exists a constant $C>2$ such that $\Phi(2x)\leq C\cdot\Phi(x)$ for all $x\in\IR$,
  abbreviated by $\Phi\in\dt$. We say that $\Phi$ obeys the global {\em $\nt$-condition}
  if the complementary $\cN$-function $\Psi$ obeys the global $\dt$-condition,
  abbreviated by $\Phi\in\nt$. Note that $\Phi\in\nt$ if and only if there exists
  a constant $c>1$ such that $\Phi(x)\leq (2c)^{-1}\Phi(cx)$ for all $x\in\IR$
  \cite[Theorem 1.1.2, p.3]{rr:02:aoos}.
\end{definition}

\begin{example}
  For $p\in(1,\infty)$ the function $\Phi_p:\IR\to\IR_+$ defined by $\Phi_p(x):=\betrag{x}^p/p$
  is an $\cN$-function and $\Phi_p\in\dt\cap\nt$. Moreover, the complementary $\cN$-function
  to $\Phi_p$ is $\Phi_q$ where $q\in(1,\infty)$ is given by $1/p+1/q=1$. $q$ is called the
  {\em conjugate index} to $p$.
\end{example}

\begin{definition}[Orlicz-Space]
  Let $\Phi$ be an $\cN$-function and $(X,\Sigma,\mu)$ be a measure space.
  Then the Orlicz space $L^\Phi(X)=L^\Phi(X,\Sigma,\mu)$ is given by
  \[ L^\Phi(X) := \klg{\su\in L^0(X):\rho_\Phi(\su/\alpha)<\infty\text{ for some }\alpha>0}
  \]
  where $\rho_\Phi$ is the Luxemburg functional given by $\rho_\Phi(\su) := \int_X \Phi(\su)\;d\mu$.
  The space $L^\Phi(X)$ endowed with the Luxemburg norm $\norm{\cdot}_\Phi$ is a Banach space
  \cite[Theorem 3.3.10, p.67]{rr:91:toos}. If in addition $\Phi\in\Delta_2\cap\nt$,
  then $L^\Phi(X)$ is reflexive \cite[Theorem 4.2.10, p.112]{rr:91:toos}.
  For reflexivity/p-convexity of Musielak-Orlicz spaces see \cite{hudzik:00:npmo}.
  Moreover, we have that $\norm{fg}_{L^1}\leq 2\norm{f}_\Phi\norm{g}_\Psi$ where
  $\Phi$ and $\Psi$ are complementary $\cN$-functions.
\end{definition}

\begin{example}[Orlicz-Sobolev spaces]\label{ex:orli-one}
  For an $\cN$-function $\Phi$ and an open set $\Omega\subset\IR^N$ we let
  $W^{1,\Phi}(\Omega)\subset L^\Phi(\Omega)\subset L^1_{\loc}(\Omega)$ be the first
  order {\em Orlicz-Sobolev space} consisting of all functions $\su$ in the Orlicz space
  $L^\Phi(\Omega)$ whose distributional derivatives of order one belong to $L^\Phi(\Omega)$.
  Then the Orlicz-Sobolev space $W^{1,\Phi}(\Omega)$ is a Banach space for the
  Luxemburg norm $\norm{\cdot}_{1,\Phi}$ associated to the Luxemburg functional
  \[ \rho_{1,\Phi,\Omega}(\su) := \int_X \Phi(\su)+\Phi(\betrag{\nabla\su})\;d\lambda
  \]
  If in addition $\Phi\in\dt\cap\nt$, then $W^{1,\Phi}(\Omega)$
  is a reflexive Luxemburg space of type $\Lambda$:
  \begin{itemize}
    \item For property \prtcc\ see Rao and Ren 
          \cite[Theorem 3.4.12, p.52 and Corollary 3.4.15]{rr:91:toos};
    \item For property \prtnr\ see Rao and Ren \cite[Corollary 5.3.4(iii), p.174]{rr:91:toos};
    \item For property \prtrn\ see Rao and Ren \cite[Theorem 1.2.7(iii), p.16]{rr:02:aoos};
    \item For reflexivity: First check that $\norm{\cdot}$ given by 
                  $\norm{\su}:=\sum\limits_{\betrag{\alpha}\leq 1}\norm{D^\alpha\su}_\Phi$
                  is equivalent to $\norm{\cdot}_{1,\Phi}$ and then
            identify $(W^{1,\Phi}(\Omega),\norm{\cdot})$ 
                  with a closed subspace of $(L^\Phi(\Omega))^{N+1}$.
  \end{itemize}
  Let $V:=\tW^{1,\Phi}(\Omega)$ be the closure of $C_c(\overline\Omega)\cap W^{1,\Phi}(\Omega)$
  in $W^{1,\Phi}(\Omega)$. Then $(V,\rho_{1,\Phi,\Omega})$ is of class $\Upsilon(\overline\Omega)$
  where the tms $X:=\overline\Omega$ is as in Example \ref{ex:tms}(2).
  Moreover, if $\Omega_1\subset\Omega_2$ are non-empty open set in $\IR^N$, then
  $\tW^{1,\Phi}(\Omega_1)\preceq\tW^{1,\Phi}(\Omega_2)$.
\end{example}

\begin{definition}[$V$-Admissibility]
  Let $(V,\rho)$ be of class $\Upsilon$. We call a continuous and bijective function
  $\psi:\IR_+\to\IR_+$ {\em $V$-admissible} if 
  \begin{itemize}
    \item[(A1)] For all $C>0$ there exists $K(C)>0$ such that
                $\psi(Ca)\leq K(C)\psi(a)$ for all $a>0$;
    \item[(A2)] For all $\su,\sv\in V$ it holds true that
                $(\psi\circ\rho)(\su\vee\sv) \leq (\psi\circ\rho)(\su)+(\psi\circ\rho)(\sv)$.
  \end{itemize}
  We will call a $V$-admissible function $\psi$ {\em strongly $V$-admissible} if for all
  $\su,\sv\in V$
  \[ (\psi\circ\rho)(\su\vee\sv)+(\psi\circ\rho)(\su\wedge\sv)\leq 
     (\psi\circ\rho)(\su)+(\psi\circ\rho)(\sv)
  \]
\end{definition}

\begin{example}\label{ex:sobo-two}
  With the assumptions and notations of Example \ref{ex:sobo-one} and $\psi(x):=x^p$ we get that
  $\psi$ is strongly $\tW^{1,p}(\Omega)$-admissible. In fact, let 
  $\su_1,\su_2\in\tW^{1,p}(\Omega)$.
  By considering the disjoint sets $D_1:=\klg{x\in\Omega:\su_1(x)<\su_2(x)}$,
  $D_2:=\klg{x\in\Omega:\su_1(x)>\su_2(x)}$ and
  $D_3:=\klg{x\in\Omega:\su_1(x)=\su_2(x)}$ we get from Stampacchia's Lemma
  \begin{eqnarray*}
    \norm{\su_1\vee \su_2}_{W^{1,p}(\Omega)}^p &=& 
        \int_{D_1} \betrag{\su_2}^p+\betrag{\nabla \su_2}^p +
        \int_{D_2} \betrag{\su_1}^p+\betrag{\nabla \su_1}^p +
        \int_{D_3} \betrag{\su_1}^p+\betrag{\nabla \su_1}^p\\
    \norm{\su_1\wedge\su_2}_{W^{1,p}(\Omega)}^p &=& 
        \int_{D_1} \betrag{\su_1}^p+\betrag{\nabla \su_1}^p +
        \int_{D_2} \betrag{\su_2}^p+\betrag{\nabla \su_2}^p +
        \int_{D_3} \betrag{\su_2}^p+\betrag{\nabla \su_2}^p.
  \end{eqnarray*}
  From this we deduce that
  \[ \psi\klrk{\norm{\su_1\vee\su_2}_{W^{1,p}(\Omega)}} + 
     \psi\klrk{\norm{\su_1\wedge\su_2}_{W^{1,p}(\Omega)}} 
      = \psi\klrk{\norm{\su_1}_{W^{1,p}(\Omega)}} + \psi\klrk{\norm{\su_2}_{W^{1,p}(\Omega)}}.
  \]
\end{example}

\begin{example}\label{ex:orli-two}
  With the assumptions and notations from Example \ref{ex:orli-one} and $\psi(x):=x$ we get that
  $\psi$ is strongly $\tW^{1,\Phi}(\Omega)$-admissible. In fact, let
  $\su_1,\su_2\in\tW^{1,\Phi}(\Omega)$. By considering the disjoint sets
  $D_1:=\klg{x\in\Omega:\su_1(x)<\su_2(x)}$,
  $D_2:=\klg{x\in\Omega:\su_1(x)>\su_2(x)}$ and
  $D_3:=\klg{x\in\Omega:\su_1(x)=\su_2(x)}$ we get from Stampacchia's Lemma 
  \[ \rho_{1,\Phi}(\su_1\vee \su_2) =
        \int_{D_1} \Phi(\su_2)+\Phi(\betrag{\nabla \su_2}) +
        \int_{D_2} \Phi(\su_1)+\Phi(\betrag{\nabla \su_1}) +
        \int_{D_3} \Phi(\su_1)+\Phi(\betrag{\nabla \su_1})
  \]
  \[
    \rho_{1,\Phi}(\su_1\wedge\su_2) = 
        \int_{D_1} \Phi(\su_1)+\Phi(\betrag{\nabla \su_1}) +
        \int_{D_2} \Phi(\su_2)+\Phi(\betrag{\nabla \su_2}) +
        \int_{D_3} \Phi(\su_2)+\Phi(\betrag{\nabla \su_2}). 
  \]
  From this we deduce that $\rho_{1,\Phi,\Omega}(\su_1\vee\su_2)+\rho_{1,\Phi,\Omega}(\su_1\wedge\su_2) =
     \rho_{1,\Phi,\Omega}(\su_1)+\rho_{1,\Phi,\Omega}(\su_2)$.
\end{example}

\begin{definition}[Cutoff-Property]
  Let $(V,\rho)$ be of class $\Upsilon$,
  $K\subset X$ compact and $U\subset X$ open containing $K$.
  A function $\eta\in V\cap C_c(U)$ is called a {\em $(K,U)$-cutoff function} if
  \begin{itemize}
    \item[\prco] $\eta\equiv 1$ on $K$ and $0\leq\eta\leq 1$ on $X$;
    \item[\prcc] The mapping $V\to V$, $\su\mapsto\eta\su$ is well-defined and continuous;
    \item[\prcr] There exists a constant $C>0$ such that $\rho(\eta\su)\leq C\cdot\rho(\su)$
                for all $\su\in V$.
  \end{itemize}
  We say that a space $V$ of class $\Upsilon$ satisfies the {\em cutoff-property} if
  for every compact set $K$ and for every open set $U$ containing $K$ 
  there exists a $(K,U)$-cutoff function.
\end{definition}

\begin{example}
  Let $\Omega\subset\IR^N$, $X:=\overline\Omega$, $U$ an open set in the tms $X$, and
  $K\subset U$ compact. Then there exists an open set $O$ in $\IR^N$ such that
  $U=O\cap X$ and a test function $\varphi\in\cD(\IR^N)=C^\infty_c(\IR^N)$ such that
  $\varphi\equiv 1$ on $K$, $0\leq \varphi\leq 1$ on $\IR^N$ and $\supp(\varphi)\subset O$.
  We remark that if $\su\in W^{1,1}_{\loc}(\Omega)$ and $\varphi\in W^{1,\infty}(\Omega)$, then
  $\su\varphi\in W^{1,1}_{\loc}(\Omega)$ and $D_j(\su\varphi)=\varphi D_j\su + \su D_j\varphi$
  in $\cD(\Omega)'$. This show that
  $\betragk{\varphi\su} \leq \normk{\varphi}_{L^\infty(\Omega)}\betragk{\su}$ and
  \[ \betragk{D_j(\varphi\su)} \leq 
        \norm{\varphi}_{L^\infty(\Omega)}\betrag{D_j\su}+
        \normk{D_j\varphi}_{L^\infty(\Omega)}\betragk{\su} \Rightarrow
     \betrag{\nabla(\varphi\su)}\leq C_1\klr{\betrag{\nabla u}+\betrag{u}}
  \]
  for some constant $C_1=C_1(\varphi,N)>0$. Hence, using that $\Phi\in\dt$, we get
  \[ \int_\Omega \Phi(\varphi\su)+\Phi(\betrag{\nabla\varphi\su})\;dx \leq
     C_2\int_\Omega \Phi(\su)+\Phi(\betrag{\nabla u}).
  \]
  This shows that for $\Phi\in\dt\cap\nt$ we get that
  $\tilde W^{1,\Phi}(\Omega)\in \Upsilon(\overline\Omega)$ has the
  cutoff-property. Similarly, we get that for $p\in(1,\infty)$ the Sobolev
  $W^{1,p}(\Omega)\in\Upsilon(\overline\Omega)$ has the cutoff-property.
\end{example}

\subsection{The $\Upsilon$-Capacity}\label{ss:rcap}

\begin{definition}[Choquet Capacity]
  Let $(X,\tau)$ be a topological space and let $\cK$ denote the collection of all
  compact sets in $X$. Then a mapping $\sC$ from the power set $\cP(X)$ of $X$
  into $[-\infty,\infty]$ is called a {\em Choquet capacity} on the paved space
  $(X,\cK)$ if
  \begin{itemize}
    \item[\prchm] If $A_1\subset A_2\subset X$, then $\sC(A_1)\leq\sC(A_2)$;
    \item[\prchi] If $(A_n)_n\subset X$ is increasing and $A=\bigcup_n A_n$, then $\lim_n \sC(A_n)=\sC(A)$;
    \item[\prchd] $K_n\in\cK$ is decreasing and $K=\bigcap_n K_n$, then $\lim_n \sC(K_n)=\sC(K)$.
  \end{itemize}
  We will call a Choquet capacity $\sC$ normed, if $\sC(\emptyset)=0$.
\end{definition}

\begin{definition}\label{def:upcap}
  Let $V$ be of class $\Upsilon$ and $\psi$ be a $V$-admissible function. Then we define the
  {\em $\Upsilon$-capacity} $\Cap_{\psi,V}$ of an arbitrary set $A\subset X$ by
  \[ \Cap_{\psi,V}(A):= \inf\klg{(\psi\circ\rho)(\su):\su\in\cY_V(A)}
  \]
  where $\cY_V(A):=\klg{\su\in V:\exists O\text{ open in }X, 
         A\subset O,\su\geq 1\text{ $\mu$-a.e. on }O}$.
  When $\Phi\in\dt\cap\nt$ is an $\cN$-function, then we call the capacity
  $\Cap_{\Phi,\Omega}:=\Cap_{\psi,V}$ the {\em relative $\Phi$-capacity}
  where $\psi(x):=x$ and $V:=\tW^{1,\Phi}(\Omega)$. When $p\in(1,\infty)$ then we call
  the capacity $\Cap_{p,\Omega}:=\Cap_{\psi,V}$ the {\em relative $p$-capacity}
  where $\psi(x):=\betrag{x}^p$ and $V:=\tW^{1,p}(\Omega)$. Note that
  when $\Omega=\IR^N$, then $\Cap_p:=\Cap_{p,\Omega}=\Cap_{p,\IR^N}$ is the classical
  $p$-capacity.
\end{definition}


\section{Properties of the $\Upsilon$-Capacity}

\subsection{Elementary Properties}

In this subsection we assume that $(V,\rho)$ is of class $\Upsilon$ and
$\psi$ is $V$-admissible.

\begin{lemma}\label{lem:weaklyclosed}
  For every open set $O\subset X$ the set $\cY_V(O)$ is convex and (weakly) closed in $V$.
\end{lemma}

\begin{proof}
  Since $\cY_V(O)=\klg{\su\in V:\su\geq 1\text{ $\mu$-a.e. on }O}$
  convexity is clear. That $\cY_V(O)$ is closed follows from property \prss.
\end{proof}

\begin{proposition}\label{prop:open-extr}
  Let $O\subset X$ be an open set with $\Cap_{\psi,V}(O)<\infty$.
  Then there exists a function $\su\in\cY_V(O)$ such that
  $\Cap_{\psi,V}(O)=\psi(\rho(\su))$, $\su=1$ $\mu$-a.e. on $O$ and
  $0\leq \su\leq 1$ $\mu$-a.e. on $X$. If in addition $\rho$ is strictly
  convex, then $\su$ is unique.
\end{proposition}

\begin{proof}
  Let $\su_n\in V$ be such that $\psi(\rho(\su_n))\to\Cap_{\psi,V}(O)=:\psi(c)$.
  Then $(\su_n)_n$ is a bounded sequence in $V$ [property \prtnr]. Since $V$ is reflexive, by possibly
  passing to a subsequence, we may assume that $(\su_n)_n$ converges weakly to a
  function $\sv\in V$. Using Lemma \ref{lem:weaklyclosed} and Remark \ref{rem:semicont}
  we get that $\sv\in\cY_V(O)$ and $\Cap_{\psi,V}(O)=\psi(\rho(\sv))$.
  Now let $\su:=(\sv\wedge 1)^+$. Then $\su\in\cY_V(O)$ and $\rho(\su)\leq\rho(\sv)$ and
  hence $\rho(\su)=\rho(\sv)$ and therefore $\Cap_{\psi,V}(O)=\psi(\rho(\su))$.
  If $\rho$ is strictly convex we get, using that $\psi$ is strictly increasing, uniqueness
  of the minimizer $\su$.
\end{proof}

\begin{theorem}\label{thm:choquet}
  The $\Upsilon$-capacity $\Cap_{\psi,V}$ is a normed Choquet capacity
  on $X$ and for every $A\subset X$ we have that
  \begin{equation}\label{eq:byopensets}
     \Cap_{\psi,V}(A) =
     \inf\klgk{\Cap_{\psi,V}(O):O\subset X\text{ open and }A\subset O}.
  \end{equation}
\end{theorem}

\begin{proof}
  Equation \eqref{eq:byopensets}, $\Cap_{\psi,V}(\emptyset)=0$ and
  $A\subset B\subset X\Rightarrow \Cap_{\psi,V}(A)\leq\Cap_{\psi,V}(B)$
  are direct consequences of Definition \ref{def:upcap}.

  Now let $(K_n)_n$ be a decreasing sequence of compact subsets in $X$ and denote by
  $K$ the intersection of all $K_n$. If $O\subset X$ is open and contains $K$, then there
  exists $n_0\in\IN$ such that $K_n\subset O$ for all $n\geq n_0$
  [property \prmh].
  Hence $\Cap_{\psi,V}(K)\leq\lim_n\Cap_{\psi,V}(K_n)\leq\Cap_{\psi,V}(O)$. Taking the infimum
  over all open sets $O$ in $X$ containing $K$ we get from Equation \eqref{eq:byopensets}
  that $\Cap_{p,\Omega}(K)=\lim_n \Cap_{p,\Omega}(K_n)$.

  Now let $(A_n)_n$ be an increasing sequence of subsets of $X$ and denote by $A$
  the union of all $A_n$. Let $s:=\lim_n \Cap_{\psi,V}(A_n)\leq\Cap_{\psi,V}(A)\in[0,\infty]$.
  To get the converse inequality we may assume that $s<\infty$ we and let $\su_n\in\cY_{V}(A_n)$
  be such that $(\psi\circ\rho)(\su_n)\leq \Cap_{\psi,V}(A_n)+1/n$.
  Therefore $(\su_n)_n$ is a bounded sequence in the reflexive Banach space $V$ and hence
  has a weakly convergent subsequence. Let $\su\in V$ denote the weak limit of this subsequence.
  By Mazur's lemma there is a sequence $(\sv_j)_j$ consisting of convex combinations
  of the $\su_n$ with $n\geq j$ which converges strongly to $\su$. By the convexity of 
  $\rho$ we get that
  \[ \rho(\sv_j) \leq \sup_{n\geq j}\rho(\su_n)\leq \psi^{-1}(s)+1/j.
  \]
  Since $\su_n\geq 1$ $\mu$-a.e.\ on $U_n$ for an open set $U_n$ containing $A_n$ we get that
  there exists an open set $W_n$ (the finite intersection of $U_j$ with $j\geq n$) containing $A_n$
  such that $\sv_n\geq 1$ $\mu$-a.e.\ on $W_n$. Since $(\sv_j)_j$ converges to $\su$ we may assume,
  by possibly passing to a subsequence, that $\normka{\betragk{\sv_{j+1}-\sv_{j}}}_V\leq 2^{-j}$. Let
  \[ \sw_j:=\sv_j+\sum_{i=j}^\infty \betrag{\sv_{i+1}-\sv_i} \geq 
     \sv_j+\sum_{i=j}^{k-1} (\sv_{i+1}-\sv_i)=\sv_k
     \quad\text{ for } k\geq j.
  \]
  Then $\sw_j\in V$ and $\sw_j\geq 1$ $\mu$-a.e.\ on $O_j$ where the open set $O_j$ is given
  by $O_j:=\bigcup_{i=j}^\infty W_i\supset A$. Since $\sw_j\to\su$ and $\sv_j\to\su$ in $V$ we get 
  using property \prtcc
  \[ \Cap_{\psi,V}(A)\leq
     \lim_j\psi(\rho(\sw_j))=\psi(\rho(\su))=\lim_j \psi(\rho(\sv_j)) \leq s=\lim_n\Cap_{\psi,V}(A_n).
  \]
\end{proof}

\begin{lemma}\label{lem:lebesgue}
  If $A\subset X$ is $\Cap_{\psi,V}$-polar (Definition \ref{def:polar}), that is, $\Cap_{\psi,V}(A)=0$, then $\mu(A)=0$.
\end{lemma}

\begin{proof}
  Let $\su_n\in\cY_V(A)$ be such that $\rho(\su_n)\to 0$.
  Then by property \prtrn $\su_n\to 0$ in $V$.
  By possibly passing to a subsequence [property \prss] we may assume that $\su_n\to 0$
  $\mu$-a.e. Since $\su_n\geq 1$ $\mu$-a.e.\ on $A$ we get that $\mu(A)=0$.
\end{proof}

\begin{proposition}\label{prop:cont}
  Assume that $V$ has the cutoff-property. Then for every compact set $K\subset X$ we have that
  \begin{eqnarray*}
      \Cap_{\psi,V}(K) 
        &=& \inf\klg{\psi(\rho(u)):u\in V\cap C_c(X), u\geq 1\text{ on }K}\\
        &=& \inf\klg{\psi(\rho(u)):u\in V\cap C(X), u\geq 1\text{ on }K}.
  \end{eqnarray*}
\end{proposition}

\begin{proof}
  Without loss of generality we may assume that $\cY_V(K)\not=\emptyset$.
  Let $\su\in\cY_V(K)$ be fixed. Then there exists an open set $U$ containing $K$
  such that $\sv:=(\su\wedge 1)^+=1$ $\mu$-a.e.\ on $U$. Let $\eta$ be a $(K,U)$-cutoff function
  and let $(v_n)_n$ be a sequence in $V\cap C_c(X)$ which converges to $\sv$ in $V$. Then
  $u_n:=\eta+(1-\eta)v_n^+$ converges in $V$ to $\eta+(1-\eta)\sv=\eta\sv+(1-\eta)\sv=\sv$
  [properties \prcc+\prnz].
  Using that $u_n\in V\cap C_c(X)$, $u_n\geq 1$ on $K$ and $\rho(\sv)\leq\rho(\su)$ [property \prn]
  we get that
  \begin{eqnarray*}
     \Cap_{\psi,V}(K) 
       &\geq& \inf\klg{\psi(\rho(u)): u\in V\cap C_c(X), u\geq 1\text{ on }K}\\
       &\geq& \inf\klg{\psi(\rho(u)): u\in V\cap C(X), u\geq 1\text{ on }K}.
  \end{eqnarray*}
  For the converse inequality we fix a function $u\in V\cap C(X)$
  such that $u\geq 1$ on $K$. Then $u_n:=(1+1/n)u\in\cY_V(K)$ and hence
  $\Cap_{\psi,V}(K) \leq \psi(\rho(u_n))\to\psi(\rho(u))$ [property \prtcc].
\end{proof}

\begin{theorem}\label{thm:strongsub}
  If $\psi$ is strongly $V$-admissible, then the $\Upsilon$-capacity is strongly subadditive,
  that is, for all $M_1,M_2\subset X$
  \begin{equation*}
     \Cap_{\psi,V}(M_1\cup M_2)+\Cap_{\psi,V}(M_1\cap M_2)\leq\Cap_{\psi,V}(M_1)+\Cap_{\psi,V}(M_2)
  \end{equation*}
\end{theorem}

\begin{proof}
  Let $\su_j\in\cY_V(M_j)$ for $j=1,2$ and let $\su:=\su_1\vee\su_2$,
  $\sv:=\su_1\wedge\su_2$. Then we have that $\su\in\cY_V(M_1\cup M_2)$, $\sv\in\cY_V(M_1\cap M_2)$ and
  \begin{eqnarray*}
     \Cap_{\psi,V}(M_1\cup M_2)+\Cap_{\psi,V}(M_1\cap M_2)
     &\leq& \psi(\rho(\su_1\vee\su_2))+\psi(\rho(\su_1\wedge\su_2))\\
     &\leq& \psi(\rho(\su_1))+\psi(\rho(\su_2)).
  \end{eqnarray*}
  Taking the infimum over all $\su_j\in\cY_V(M_1)$ the claim follows.
\end{proof}

\begin{theorem}\label{thm:subadd}
  The $\Upsilon$-capacity is subadditive, that is, for all $M_1,M_2\subset X$
  \begin{equation*}
     \Cap_{\psi,V}(M_1\cup M_2)\leq\Cap_{\psi,V}(M_1)+\Cap_{\psi,V}(M_2)
  \end{equation*}
\end{theorem}

\begin{proof}
  The proof follows the lines in the proof of Theorem \ref{thm:strongsub}. Note that the
  $\Upsilon$-capacity $\Cap_{\psi,V}$ is defined only for $V$-admissible $\psi$ and
  that the $V$-admissibility of $\psi$ was assumed at the beginning of this section.
\end{proof}

\begin{theorem}\label{thm:csubadd}
  The $\Upsilon$-capacity $\Cap_{\psi,V}$ is countably subadditive, that is, for all $A_k\subset X$
  \[ \Cap_{\psi,V}\klr{\bigcup\nolimits_{k\in\IN} A_k} \leq \sum_{k\in\IN} \Cap_{\psi,V}(A_k).
  \]
\end{theorem}

\begin{proof}
  Let $B_n$ be the union of $A_k$ with $1\leq k\leq n$ and let $A$ be the union of all $A_k$.
  From Theorem \ref{thm:subadd} we get by induction that for all $n\in\IN$
  $\Cap_{\psi,V}\klr{B_n} \leq \sum_{k=1}^n \Cap_{\psi,V}(A_k)$.
  Using that $\Cap_{\psi,V}$ is a Choquet capacity [property \prchi] we get
  \[ \Cap_{\psi,V}\klr{A} = \lim_n \Cap_{\psi,V}\klr{B_n} 
     \leq \lim_n \sum_{k=1}^n \Cap_{\psi,V}(A_k) = \sum_{k\in\IN} \Cap_{\psi,V}(A_k).
  \]
\end{proof}

\subsection{Relations between $\Upsilon$-Capacities}

In this subsection we assume that for $j=1,2$ the tms $X_j=(X_j,\tau_j,\Sigma_j,\mu_j)$
is of type $\Theta$, $X_1\preceq X_2$, $(V_j,\rho_j)$ is of class $\Upsilon(X_j)$, $\psi_j$ is
$V_j$-admissible. Moreover, for $1\leq i,j\leq 2$ we let $\psi_{i,j}:\IR_+\to\IR_+$ be the
bijective and continuous function given by $\psi_{i,j}:=\psi_i\circ\psi_j^{-1}$.

\begin{definition}[$(V_1,V_2)$-extension property]
  We will say that a set $U\subset X_1$ has the {\em $(V_1,V_2)$-extension property} if
  \begin{itemize}
    \item[\prxm] $\mu_1(N)=0\Rightarrow \mu_2(N)=0$ for all $\Sigma_2$-measurable sets $N\subset U$;
    \item[\prxe] There exists a constant $c>0$ and a mapping $\cE:V_1^0(U)\to V_2$ where $V_1^0(U)$
                is given by
                \[ V_1^0(U):=\klg{\su\in V_1:\su=0\; \mu_1\text{-a.e. on }X_1\setminus U},
                \]
                such that $\cE\su=\su$ $\mu_1$-a.e.\ on $X_1$ and $\rho_2(\cE\su)\leq c\rho_1(\su)$.
  \end{itemize}
  We say that a set $U$ in $X_1$ has the {\em $\sigma(V_1,V_2)$-extension property} if
  $U$ is the countable union of open sets in $X_1$ satisfying the $(V_1,V_2)$-extension property.
\end{definition}

\begin{example}\label{ex:extprop}
  Let $\Omega_1\subset\Omega_2$ be non-empty open sets in $\IR^N$ and $U\subset\subset\Omega_1$
  (this means $U$ is open, $\overline U$ is compact and $\overline U\subset\Omega_1$).
  \begin{enumerate}
    \item If $1<q\leq p<\infty$, then $U$ has the
          $(\tW^{1,p}(\Omega_1),\tW^{1,q}(\Omega_2))$-extension property and hence
          $\Omega_1$ has the $\sigma(\tW^{1,p}(\Omega_1),\tW^{1,q}(\Omega_2))$-extension property.
    \item If $\Phi\in\dt\cap\nt$ is an $\cN$-function, then $U$ has the
          $(\tW^{1,\Phi}(\Omega_1),\tW^{1,\Phi}(\Omega_2))$-extension property and hence
          $\Omega_1$ has the $\sigma(\tW^{1,\Phi}(\Omega_1),\tW^{1,\Phi}(\Omega_2))$-extension property.
  \end{enumerate}
\end{example}

\begin{lemma}\label{lem:compacts}
  Assume that $U$ is open in $X_1$ and has the $(V_1,V_2)$-extension property,
  $K\subset U$ is compact and $V_1$ has the cutoff property.
  Then there exists a constant $C=C(K)>0$ such that for all $A\subset K$
  \begin{equation*}
     \Cap_{\psi_2,V_2}(A) \leq \psi_{2,1}(C\cdot\Cap_{\psi_1,V_1}(A)).
  \end{equation*}
  In particular, when $\psi:=\psi_1=\psi_2$, then there exists a constant $C>0$ such that
  \begin{equation*}
     \Cap_{\psi,V_2}(A) \leq C\cdot\Cap_{\psi,V_1}(A).
  \end{equation*}
\end{lemma}

\begin{proof}
  Let $\eta\in V_1\cap C_c(U)$ be a $(K,U)$-cutoff function and $\su\in\cY_{V_1}(A)$ be fixed.
  Then $\eta\su\in V_1^0(U)$, $\sv:=\cE(\eta\su)\in\cY_{V_2}(A)$ and hence
  [properties \prcr+\prxe]
  \begin{eqnarray*} 
     \Cap_{\psi_2,V_2}(A)
        &\leq& (\psi_2\circ\rho_2)(\sv) \leq \psi_2(c_1\rho_1(\eta\su))
         \leq \psi_2(c_2\rho_1(\su)) \\
        &=& \psi_{2,1}\klr{\psi_1(c_2\rho_1(\su))} 
         \leq \psi_{2,1}(C\cdot(\psi_1\circ\rho_1)(\su)).
  \end{eqnarray*}
  Taking the infimum over all $\su\in\cY_{\psi,V_1}(A)$ we get the claim.
\end{proof}

\begin{example}
  Let $\Omega_1\subset\Omega_2$ be non-empty open sets in $\IR^N$,
  $1<q\leq p<\infty$ and let $K\subset\Omega_1$ be compact. Then there
  exists a constant $C>0$ such that for all $A\subset K$
  \[ \Cap_{q,\Omega_2}(A) \leq C\cdot\Cap_{p,\Omega_1}(A)^{q/p}.
  \]
  In fact, let $U\subset\subset\Omega_1$ be such that $K\subset U$. Then
  by Example \ref{ex:extprop}(2) we get that all assumptions from 
  Lemma \ref{lem:compacts} are satisfied.
\end{example}

\begin{lemma}\label{lem:capesti}
  If $V_1\preceq V_2$ then there exists a $C>0$ such that for all $A\subset X_1$
  \begin{equation}\label{eq:domi}
     \Cap_{\psi_1,V_1}(A) \leq \psi_{1,2}(C\cdot\Cap_{\psi_2,V_2}(A)).
  \end{equation}
\end{lemma}

\begin{proof}
  Let $\su\in\cY_{V_2}(A)$. Then $\su|_{X_1}\in\cY_{V_1}(A)$ and
  $\rho_1(\su|_{X_1})\leq c\rho_2(\su)$ by Definition \ref{def:domups}. Hence
  \begin{eqnarray*}
     \Cap_{\psi_1,V_1}(A)
        &\leq& (\psi_1\circ\rho_1)(\su|_{X_1})\leq \psi_1(c\cdot\rho_2(\su)) \\
        &=& \psi_{1,2}(\psi_2(c\cdot\rho_2(\su))) \leq \psi_{1,2}(C\cdot(\psi_2\circ\rho_2)(\su)).
  \end{eqnarray*}
  Taking the infimum over all $\su\in\cY_{V_2}(A)$ we get the claim.
\end{proof}

\begin{proposition}\label{prop:polars}
  Assume that $U$ is an open set in $X_1$ and has the $\sigma(V_1,V_2)$-extension property,
  $V_1$ has the cutoff property and $V_1\preceq V_2$. Then for every $A\subset U$ we have that
  \begin{equation}\label{eq:polars}
     \Cap_{\psi_1,V_1}(A)=0\quad\Longleftrightarrow\quad \Cap_{\psi_2,V_2}(A)=0.
  \end{equation}
\end{proposition}

\begin{proof}
  From Equation \eqref{eq:domi} we get that $\Cap_{\psi_2,V_2}(A)=0$ implies that
  $\Cap_{\psi_1,V_1}(A)=0$. Hence to prove \eqref{eq:polars} it remains to prove the
  converse implication.
  For this let $U_m\subset X_1$ be open sets with the $(V_1,V_2)$-extension property
  such that $U=\bigcup_m U_m$ and let $K_{m,n}\subset U_m$ be compact sets
  such that $U_m=\bigcup_n K_{m,n}$ and assume that $\Cap_{\psi_1,V_1}(A)=0$.
  Then by Lemma \ref{lem:compacts} there exist constants $C_{n,m}$ such that
  \[ \Cap_{\psi_2,V_2}(K_{m,n}\cap A)\leq \psi_{2,1}(C_{m,n}\cdot\Cap_{\psi_1,V_1}(K_{m,n}\cap A))=0.
  \]
  Using that $\Cap_{\psi_2,V_2}$ is a countably subadditive [Theorem \ref{thm:csubadd}]
  and that $\bigcup_{n,m} K_{n,m}=U$ we get that
  \[ \Cap_{\psi_2,V_2}(A)\leq \sum_{m,n} \Cap_{\psi_2,V_2}(K_{n,m}\cap A) =0.
  \]
\end{proof}

\begin{corollary}\label{cor:polars}
  Let $\Omega_1\subset\Omega_2$ be non-empty open sets in $\IR^N$,
  $p\in(1,\infty)$ and $\Phi\in\dt\cap\nt$ be an $\cN$-function.
  Then for all sets in $A\subset\Omega_1$ we have that
  \[ \Cap_{p,\Omega_1}(A)=0\Leftrightarrow \Cap_{p,\Omega_2}(A)=0\quad\text{ and }\quad
     \Cap_{\Phi,\Omega_1}(A)=0\Leftrightarrow\Cap_{\Phi,\Omega_2}(A)=0.
  \]
\end{corollary}

\begin{remark}
  In general the assertion of Corollary \ref{cor:polars} (and hence of Proposition \ref{prop:polars})
  is untrue for $A\subset\overline\Omega_1$.
  To see this we let $N\geq 2$ and
  $\Omega:=(0,1)^N\setminus\bigcup_{j\in\IN} [2^{-2j},2^{1-2j}]\times[2^{-j},1]\times(0,1)^{N-2}$.
  Then the open set $\Omega\subset\IR^N$ is bounded and connected
  and $A:=\klg{0}\times[0,1]^{n-1}\subset\partial\Omega$ is such that the $(N-1)$-dimensional Hausdorff
  measure $\mathcal{H}^{N-1}(A)=1$.
  Then for every $p\in(1,\infty)$ we have that $\Cap_p(A)>0$ and $\Cap_{p,\Omega}(A)=0$.
\end{remark}

\begin{definition}[Continuous Extension Property]\label{def:context}
  We will say that $V_1$ has the {\em continuous $V_2$-extension property} if there
  exists a (possibly non-linear) mapping $\cE:V_1\to V_2$ and a constant $C=C(\cE)$ such that
  $\cE(V_1\cap C(X_1))\subset V_2\cap C(X_2)$ and for all $\su_1\in V_1$ it holds true that
  $\rho_2(\cE(\su_1))\leq C\rho_1(\su_1)$ and $\cE\su_1=\su_1$ $\mu_1$-a.e.\ on $X_1$.
\end{definition}

\begin{definition}
  Let $\Omega_1\subset\Omega_2$ be non-empty open sets in $\IR^N$ and $p\in(1,\infty)$.
  Then we say that $W^{1,p}(\Omega_1)$ has the {\em $W^{1,p}(\Omega_2)$-extension property}
  if the restriction $W^{1,p}(\Omega_2)\to W^{1,p}(\Omega_1)$, $\su\mapsto \su|_{\Omega_1}$
  is surjective. If $\Omega_2=\IR^N$ and $\Omega_1$ has the $W^{1,p}(\Omega_2)$-extension
  property, then we say briefly that $\Omega_1$ has the {\em $W^{1,p}$-extension property}.
\end{definition}

The following is an immediate consequence of Shvartsman \cite{shvartsman:07:eos} and
Haj{\l}asz and Koskela and Tuominen \cite{koskela:08:see}.

\begin{theorem}\label{thm:koskela}
  Let $p\in(1,\infty)$ and $\Omega_1\subset\IR^N$ be a $W^{1,p}$-extension domain.
  Then $\Omega_1$ has the continuous $W^{1,p}(\Omega_2)$-extension property for every
  open set $\Omega_2\subset\IR^N$ containing $\Omega_1$.
\end{theorem}

\begin{proof}
  Since $\Omega$ is a $(1,p)$-extension domain, we get from
  \cite[Theorem 2 and Lemma 2.1]{koskela:08:see} that there exists a constant
  $\delta_\Omega>0$ such that $\lambda(B(x,r)\cap\Omega)\geq \delta_\Omega r^N$
  for all $0<r\leq 1$ and $\lambda(\partial\Omega)=0$.
  For a measurable set $A\subset\IR$ we let $M^{1,p}(A)$ be
  the Sobolev-type space introduced by Haj{\l}asz consisting of
  those function $\su\in L^p(A)$ with generalized gradient in $L^p(A)$.
  It follows from \cite[Theorem 1.3]{shvartsman:07:eos} that
  $M^{1,p}(\IR^N)|_{\overline\Omega}=M^{1,p}(\overline\Omega)$ and that
  there exists a linear continuous extension operator
  $\cE:M^{1,p}(\overline\Omega)\to M^{1,p}(\IR^N)$. Using that
  $M^{1,p}(\IR^N)=W^{1,p}(\IR^N)$ as sets with equivalent norms,
  we get that $M^{1,p}(\overline\Omega)=W^{1,p}(\Omega)$
  are equal as sets with equivalent norms and hence the extension
  operator $\cE$ constructed for $M^{1,p}(\overline\Omega)$
  is also a linear continuous extension operator
  from $W^{1,p}(\Omega)$ into $W^{1,p}(\IR^N)$. It is left to the reader
  to verify that the extension operator $\cE$ constructed by
  Shvartsman \cite[Equation (1.5)]{shvartsman:07:eos} maps
  $W^{1,p}(\Omega)\cap C(\overline\Omega)$ into $W^{1,p}(\IR^N)\cap C(\IR^N)$.
\end{proof}

\begin{theorem}
  If $V_1\preceq V_2$, $V_1$ has the continuous $V_2$-extension property and the cutoff property,
  then there exist constants $C_1,C_2>0$ such that for every set $A\subset X_1$
  \[ \Cap_{\psi_2,V_2}(A)\leq \psi_{2,1}\klek{C_1\cdot \Cap_{\psi_1,V_1}(A)} 
                         \leq \psi_{2,1}\klek{C_1\cdot \psi_{1,2}\klrk{C_2\Cap_{\psi_2,V_2}(A)}}.
  \]
\end{theorem}

\begin{proof}
  Let $K\subset X_1$ be a compact set. By Proposition \ref{prop:cont} there exist
  $u_n\in V_1\cap C_c(X_1)$ such that $u_n\geq 1$ on $K$ and
  $\psi_1(\rho_1(u_n))\to\Cap_{\psi_1,V_1}(K)$. Let $\cE:V_1\to V_2$ be a continuous
  extension operator and define $v_n:=\cE u_n$.
  Then $v_n\in V_2\cap C(X_2)$ and $v_n\geq 1$ on $K$. Hence by Proposition \ref{prop:cont}
  we get that
  \begin{eqnarray*} 
     \Cap_{\psi_2,V_2}(K) 
       &\leq& \psi_2(\rho_2(v_n)) \leq \psi_2(C(\cE)\rho_1(u_n)) = \psi_{2,1}(\psi_1(C(\cE)\rho_1(u_n)))\\
       &\leq& \psi_{2,1}(C(\cE,\psi_1) \psi_1(\rho_1(u_n)))\rightarrow 
              \psi_{2,1}(C(\cE,\psi_1)\Cap_{\psi_1,V_1}(K)).
  \end{eqnarray*}
  Let $W$ be an open set in $X_1$. Then there exists an increasing sequence $(K_n)_n$ of
  compact sets such that $\bigcup_n K_n=W$ [property \prmc]. Using that $\Cap_{\psi_1,V_1}$ and
  $\Cap_{\psi_2,V_2}$ are Choquet capacities [property \prchi] we get that
  \begin{eqnarray*}
     \Cap_{\psi_2,V_2}(W)
       &=& \lim_n \Cap_{\psi_2,V_2}(K_n) \leq 
           \lim_n \psi_{2,1}(C(\cE,\psi_1)\Cap_{\psi_1,V_1}(K_n)) \\
       &=& \psi_{2,1}(C(\cE,\psi)\Cap_{\psi,V_1}(W)).
  \end{eqnarray*}
  Now let $A\subset X_1$ be arbitrary. Then by Theorem \ref{thm:choquet}
  \begin{eqnarray*}
    \Cap_{\psi_2,V_2}(A) 
        &=& \inf\klgk{\Cap_{\psi_2,V_2}(O):O\text{ is open in }X_2\text{ and } A\subset O} \\
        &=& \inf\klgk{\Cap_{\psi_2,V_2}(O\cap X_1):O\text{ is open in }X_2\text{ and } A\subset O} \\
        &=& \inf\klgk{\Cap_{\psi_2,V_2}(W):W\text{ is open in }X_1\text{ and } A\subset W} \\
        &\leq& \psi_{2,1}\klek{C(\cE,\psi_1) \inf\klgk{\Cap_{\psi_1,V_1}(W):
                     W\text{ is open in }X_1\text{ and } A\subset W}}\\
        &=& \psi_{2,1}\klek{C(\cE,\psi)\Cap_{\psi,V_1}(A)}.
  \end{eqnarray*}
  The remaining inequality follows from Lemma \ref{lem:capesti}.
\end{proof}


\subsection{Quasicontinuity and Polar Sets}

In this subsection we assume that the tms $X=(X,\tau,\Sigma,\mu)$ is of type $\Theta$,
$(V,\rho)$ is of class $\Upsilon(X)$ and $\psi$ is $V$-admissible. The purpose of this
subsection is to prove existence and uniqueness of $\Cap_{\psi,V}$-quasi continuous
representatives on $X$.

\begin{definition}\label{def:polar}
  A set $P\subset X$ is said to be {\em $\Cap_{\psi,V}$-polar} if $\Cap_{\psi,V}(P)=0$.
  A pointwise defined function $u$ on $\domain\subset X$ is called 
  {\em $\Cap_{p,\Omega}$-quasi continuous} on $\domain$ if for each $\varepsilon>0$
  there exists an open set $O$ in $X$ with $\Cap_{\psi,V}(O)<\varepsilon$ such that
  $u$ restricted to $\domain\setminus O$ is continuous. We say that a property
  holds $\Cap_{\psi,V}$-quasi everywhere (briefly $\Cap_{\psi,V}$-q.e.) if
  it holds except for a $\Cap_{\psi,V}$-polar set.
\end{definition}

\begin{lemma}\label{lem:conver}
  If $\su\in V$ and $u_k\in V\cap C_c(X)$ are such that
  \[ \sum_{k=1}^\infty \psi(2^{k+2}\normk{\su-u_{k+1}}_V+2^{k+1}\norm{\su-u_k}_V)<\infty,
  \]
  then the pointwise limit $\tilde u:=\lim_k u_k$ exists $\Cap_{\psi,V}$-quasi everywhere
  on $X$, $\tilde u:X\to\IR$ is $\Cap_{\psi,V}$-quasi continuous 
  and $\tilde u=\su$ $\mu$-almost everywhere on $X$.
\end{lemma}

\begin{proof}
  Let $k_0\in\IN$ be such that $2^k\normk{\su-u_k}_V<1/4$ for all $k\geq k_0$ and 
  consider the open sets $G_k:=\klgk{x\in X:\betragk{u_{k+1}(x)-u_k(x)}>2^{-k}}$.
  Then $2^k\cdot\betragk{u_{k+1}-u_k}\geq 1$ on $G_k$ and
  \[ \normk{u_{k+1}-u_k}_V\leq
     \normk{\su-u_{k+1}}_V + \normk{\su-u_k}_V \leq
     2\normk{\su-u_{k+1}}_V + \normk{\su-u_k}_V.
  \]
  Therefore, if $k\geq k_0$, we get
  \[ \normka{2^k\betrag{u_{k+1}-u_k}}_V\leq \normk{2^{k+1}(u_{k+1}-u_k)}_V\leq 
        2^{k+2}\norm{\su-u_{k+1}}_V+2^{k+1}\norm{\su-u_k}_V<1.
  \]
  From this we deduce that (see Remark \ref{rem:semicont}) for all $k\geq k_0$
  \begin{eqnarray*} 
     \Cap_{\psi,V}(G_k)
       &\leq& \psi\klek{\rho(2^k\cdot\betrag{u_{k+1}-u_k})} \leq
              \psi(\normka{2^k\betrag{u_{k+1}-u_k}}_V) \\
       &\leq& \psi\klrk{ 2^{k+2}\norm{\su-u_{k+1}}_V + 2^{k+1}\norm{\su-u_k}_V}
  \end{eqnarray*}
  and hence $\sum_{k} \Cap_{\psi,V}(G_k)<\infty$. Given $\varepsilon>0$ there exists $k_1\geq k_0$
  such that $\Cap_{\psi,V}(G)<\varepsilon$ where $G:=\bigcup_{k\geq k_1} G_k$.
  Since $\betrag{u_{k+1}-u_k}\leq 2^{-k}$ on $X\setminus G$ for all $k\geq k_1$ we have
  that $(u_k)_k$ is a sequence of continuous functions on $X$ which converges
  uniformly on $X\setminus G$. Since $\varepsilon>0$ was arbitrary
  we get that $\tilde u:=\lim_k u_k$ exists $\Cap_{\psi,V}$-quasi everywhere on $X$
  and $\tilde u|_{X\setminus G}$ is continuous. To see that $\tilde u$ coincides with
  $\su$ $\mu$-almost everywhere on $X$ we argue as follows. Since $(u_k)_k$
  converges to $\su$, by possibly passing to a subsequence [property \prss] we have that
  $u_k$ converges to $\su$ $\mu$-almost everywhere. Since $(u_k)_k$ converges to
  $\tilde u$ $\Cap_{p,\Omega}$-quasi everywhere on $X$ (and hence $\mu$-almost everywhere on $X$)
  we get that $\tilde u=\su$ $\mu$-almost everywhere on $X$ (see Lemma \ref{lem:lebesgue}).
\end{proof}

\begin{theorem}\label{thm:cont-version}
  For every $\su\in V$ there exists a $\Cap_{\psi,V}$-quasi
  continuous function $\tilde u:X\to\IR$ such that
  $\tilde u=\su$ $\mu$-a.e.\ on $X$, that is, $\tilde u\in\su$.
\end{theorem}

\begin{proof}
  Let $\su\in V$. Then by definition there exists a sequence $u_n\in V\cap C_c(X)$ such
  that $u_n\to \su$ in $V$. Then a subsequence of $(u_n)_n$ satisfies the assumptions
  of Lemma \ref{lem:conver}.
\end{proof}

\begin{lemma}\label{lem:approx}
  Let $A\subset X$, $\su\in V$ be non-negative and let $u\in\su$ be a
  $\Cap_{\psi,V}$-quasi continuous version of $\su$ such that $u\geq 1$
  $\Cap_{\psi,V}$-quasi everywhere on $A$. Then there is a sequence
  $(\su_n)_n\subset \cY_V(A)$ which converges to $\su$ in $V$.
\end{lemma}

\begin{proof}
  Let $O_n$ be an open set in $X$ such that $u|_{X\setminus O_n}$
  is continuous, $u\geq 1$ everywhere on $A\setminus O_n$ and $\Cap_{\psi,V}(O)\leq 1/n$.
  Let $\sv_n$ be a capacitary extremal for $O_n$ (see Proposition \ref{prop:open-extr}) 
  such that $0\leq \sv_n\leq 1$ $\mu$-a.e.\ on $X$. Then we have that
  \[ \psi(\rho(\sv_n))=\Cap_{\psi,V}(O_n)\to 0 \Rightarrow \rho(\sv_n)\to 0 \Rightarrow \sv_n\to 0\text{ in } V.
  \]
  Let $\sw_n:=(1+1/n)u+\sv_n\geq\sv_n$. Then $\sw_n\to \su$ in $V$.
  For the open set $G_n$ in $X$ given by
  \[ G_n:=O_n\cup \klg{x\in X\setminus O_n:u(x)>n/(n+1)}
  \]
  we have that $\sw_n\geq 1$ $\mu$-a.e.\ on $G_n$ and $A\subset G_n$ and hence
  $\sw_n\in\cY_{V}(A)$.
\end{proof}

\begin{lemma}\label{lem:normest}
  Let $u\in\su\in V$ be a $\Cap_{\psi,V}$-quasi continuous version of $\su$ and let
  $a\in(0,\infty)$. Then
  \[ \Cap_{\psi,V}\klrk{\klgk{x\in X:u(x)>a}}\leq \psi\kle{\rho(\su/a)}
  \]
\end{lemma}

\begin{proof}
  Let $A:=\klgk{x\in X:u(x)>a}$. By Lemma \ref{lem:approx} there exists a sequence
  $(\su_n)_n\in\cY_V(A)$ which converges to $a^{-1}\su^+$ in $V$.
  Note that $u^+$ is a $\Cap_{\psi,V}$-quasi continuous version of $\su^+$. Hence
  \[ \Cap_{\psi,V}\klrk{\klg{x\in X:u(x)>a}}\leq
     \psi\kle{\rho(\su_n)}\to \psi\kle{\rho(\su^+/a)}\leq \psi\kle{\rho(\su/a)}.
  \]
\end{proof}

\begin{theorem}\label{thm:esti}
  Assume that $V$ has the cutoff-property and let $\su,\sv\in V$ be such that
  $\su\leq\sv$ $\mu$-a.e.\ on $U$ where $U$ is an open set in $X$.
  If $u\in\su$ and $v\in\sv$ are $\Cap_{\psi,V}$-quasi continuous
  versions of $\su$ and $\sv$, respectively,
  then $u\leq v$ $\Cap_{\psi,V}$-quasi everywhere on $U$.
\end{theorem}

\begin{proof}
  Let $(K_n)_n$ be a sequence of compact sets such that $U=\bigcup_n K_n$.
  For the sequence of compact sets we choose non-negative $(K_n,U)$-cutoff functions
  $\varphi_n\in V\cap C_c(U)$. Then the function $w_n:=\varphi_n(u-v)^+=0$ $\mu$-a.e.\ on $X$
  and we get by Lemma \ref{lem:normest}, using that $\varphi_n (u-v)^+$ is $\Cap_{\psi,V}$-quasi
  continuous, that $w_n=0$ $\Cap_{p,\Omega}$-quasi everywhere on $X$ and hence that
  $u\leq v$ $\Cap_{\psi,V}$-quasi everywhere on $K_n$ for each $n\in\IN$.
  Since the countable union of $\Cap_{\psi,V}$-polar sets is $\Cap_{\psi,V}$-polar
  we get that $u\leq v$ $\Cap_{\psi,V}$-quasi everywhere on $U$.
\end{proof}

\begin{theorem}\label{thm:quasicontrep}
  Let $\su\in V$. Then there exists a unique (up to a $\Cap_{\psi,V}$-polar set)
  $\Cap_{\psi,V}$-quasi continuous function $\tilde u:X\to\IR$ such that
  $\su=\tilde u$ $\mu$-a.e.\ on $\Omega$.
\end{theorem}

\begin{proof}
  The existence follows from Theorem \ref{thm:cont-version}. To show uniqueness we let
  $u_1,u_2\in\su\in V$ be two quasi-continuous versions of $\su$.
  Then $u_1=u_2$ $\mu$-a.e.\ on $X$ and hence by Lemma \ref{lem:normest} we get that
  $u_1=u_2$ $\Cap_{\psi,V}$-quasi everywhere on $X$.
\end{proof}

\begin{definition}
  By $\cN^\star(V,\rho)$ we denote the set of all $\Cap_{\psi,V}$-polar sets in $X$
  and we denote by $C(V,\rho)$ the space of all $\Cap_{\psi,V}$-quasi continuous functions
  $u:X\to\IR$. Note that $\cN^\star(V,\rho)$ and $C(V,\rho)$ do not depend on the
  $V$-admissible $\psi$. On $C(V,\rho)$ we define the equivalence relation $\sim$ by
  \[ u\sim v\quad :\Leftrightarrow \quad \exists P\in\cN^\star(V,\rho): u=v \text{ everywhere on }
     X\setminus P.
  \]
  For a function $u\in C(V,\rho)$ we denote by $[u]$ the equivalence class of $u$ with respect
  to $\sim$. Now the refined space $\sV$ is defined by
  $\sV := \klg{\tilde\su:\su\in V}\subset C(V,\rho)/\sim$
  where $\tilde\su:=[u]$ with $u\in\su\in V$ $\Cap_{\psi,V}$-quasi continuous.
  We equip $\sV$ with the norm $\norm{\cdot}_V$. Note that by Theorem \ref{thm:quasicontrep}
  $\sV$ is isometrically isomorphic to $V$.
  For a sequence $(\su_n)_n$ in $\sV$ and $\su\in\sV$ we say that
  {\em $(\su_n)_n$ converges $\Cap_{\psi,V}$-quasi everywhere to $\su$} if
  for every $u_n\in\su_n$ and $u\in\su$ there exists a $\Cap_{\psi,V}$-polar set $P$ such that
  $u_n\to u$ everywhere on $X\setminus P$.
  We say that {\em $(\su_n)_n$ converges $\Cap_{\psi,V}$-quasi uniformly to $\su$}
  if for every $u_n\in\su_n$, $u\in\su$ and $\varepsilon>0$ there exists an open set
  $G$ in $X$ such that $\Cap_{\psi,V}(G)\leq\varepsilon$ and $u_n\to u$ uniformly (everywhere) on
  $X\setminus G$.
\end{definition}

\begin{theorem}\label{thm:convsubse}
  If $\su_n\in\sV$ converges to $\su\in\sV$ in $\sV$, then there exists a subsequence
  which converges $\Cap_{\psi,V}$-quasi everywhere and -quasi uniformly on $X$ to $\su$.
\end{theorem}

\begin{proof}
  By possibly passing to a subsequence we may assume that
  $\sum_{n\in\IN} \psi\kle{\rho(n\betrag{\su_n-\su})}<\infty$.
  We show under the above assumption that $(\su_n)_n$ converges $\Cap_{\psi,V}$-quasi everywhere 
  and -quasi uniformly on $X$ to $\su$. Let $u_n\in\su_n$ and $u\in\su$ be fixed and define
  \[ G_n:=\klg{x\in X:\betrag{u_n(x)-u(x)}>n^{-1}}.
  \]
  We show that $u_n(x)\to u(x)$ for all $x\in X\setminus P$ where
  $P:=\bigcap_{j=1}^\infty \bigcup_{k=j}^\infty G_k$.
  If $x\in X\setminus P$ then there exists $j_0\in\IN$ such that
  $x\not\in \bigcup_{k=j_0}^\infty G_k$, that is, $\betrag{u_n(x)-u(x)}\leq n^{-1}$
  for all $n\geq j_0$ and hence $u_n\to u$ uniformly on
  $X\setminus\bigcup_{k=j_0}^\infty G_k\ni x$ and everywhere on $X\setminus P$.
  We show that $P$ is a $\Cap_{\psi,V}$-polar set.
  Let $\varepsilon>0$ be given. Then there exists $N=N(\varepsilon)$ such that
  $\sum_{n=N}^\infty \psi\kle{\rho(n\betrag{\su_n-\su})}\leq\varepsilon$. By Lemma \ref{lem:normest}
  we get that
  \[ \Cap_{\psi,V}\klr{\bigcup\nolimits_{n\geq N} G_n} \leq
     \sum_{n=N}^\infty \Cap_{\psi,V}(G_n)\leq 
     \sum_{n=N}^\infty \psi\kle{\rho(n\cdot\betrag{\su_n-\su})}\leq\varepsilon.
  \]
  Therefore $\Cap_{p,\Omega}(P)\leq\varepsilon$ and since $\varepsilon>0$
  was arbitrary the claim follows.
\end{proof}

\begin{lemma}\label{lem:loccont}
  Let $U\subset X$ be a non-empty open set which is the countable union
  of open sets $\omega_n\subset U$, $n\in\IN$. Then a function $u:X\to\IR$
  is $\Cap_{\psi,V}$-quasi continuous on $U$ if and only if $u$ is
  $\Cap_{\psi,V}$-quasi continuous on every set $\omega_n$.
\end{lemma}

\begin{proof}
  Assume that $u$ is $\Cap_{\psi,V}$-quasi continuous on every set $\omega_n\subset U$ and
  let $\varepsilon>0$ be given. Then there exists an open set $O_n\subset X$ such that
  $\Cap_{\psi,V}(O_n)\leq \varepsilon 2^{-n}$ and $u|_{\omega_n\setminus O_n}$ is continuous, that is,
  $u^{-1}(I)\cap (\omega_n\setminus O_n)$ is open in $\omega_n\setminus O_n$
  for every open set $I\subset\IR$.
  Let $O:=\bigcup_n O_n$. Then $\Cap_{\psi,V}(O)\leq \sum_n \Cap_{\psi,V}(O_n)\leq\varepsilon$ and
  $u|_{U\setminus O}$ is continuous.
  In fact, using that $u^{-1}\cap(\omega_n\setminus O)$ is open in $\omega_n\setminus O$ and hence
  open in $U\setminus O$ we get that
  $u^{-1}(I)\cap (U\setminus O) = \bigcup_n u^{-1}(I)\cap(\omega_n\setminus O)$
  is open in $U\setminus O$ for every open set $I\subset\IR$, hence
  $u|_{U\setminus O}$ is continuous.
\end{proof}

\begin{theorem}\label{thm:quasicont}
  For $j=1,2$ we let $X_j$ be a tms of type $\Theta$, $V_j$ be of class $\Upsilon(X_j)$,
  $\psi_j$ be $V_j$-admissible. If $U\subset X_1$ has the $\sigma(V_1,V_2)$-extension property,
  $V_1$ has the cutoff-property and $V_1\preceq V_2$, then a function
  $u:U\to\IR$ is $\Cap_{\psi_1,V_1}$-quasi continuous on $U$ if and only if it is
  $\Cap_{\psi_2,V_2}$-quasi continuous on $U$.
\end{theorem}

\begin{proof}
  If $u$ is $\Cap_{\psi_2,V_2}$-quasi continuous, then $u$ is $\Cap_{\psi_1,V_1}$-quasi continuous
  by Lemma \ref{lem:capesti}. Assume now that $u$ is $\Cap_{\psi_1,V_1}$-quasi continuous and let
  $K_n\subset U$ be an increasing sequence of compact sets such that $U=\bigcup_n K_n$.
  Since $V_1$ has the cutoff-property we let $\varphi_n\in V_1\cap C_c(U)$ be a
  $(K_n,U)$-cutoff function and $\omega_n:=\klg{x\in X_1:\varphi_n(x)>0}$ be an open set. Then
  $K_n\subset\omega_n\subset\subset U$. By Lemma \ref{lem:compacts} there exists
  a constant $C_n$ such that $\Cap_{\psi_2,V_2}(A)\leq \psi_{2,1}(C_n\Cap_{\psi_1,V_1}(A))$
  for all $A\subset w_n$. Hence $u$ is $\Cap_{\psi_2,V_2}$-quasi continuous on $\omega_n$.
  Lemma \ref{lem:loccont} shows now that $u$ is $\Cap_{\psi_2,V_2}$-quasi continuous on $U$.
\end{proof}

\begin{example}
  Let $\Omega\subset\IR^N$ be an open set.
  If $p\in(1,\infty)$ and $u\in\su\in\tW^{1,p}(\Omega)$ is a 
  $\Cap_{p,\Omega}$-quasi continuous version of $\su$, then $u$ is
  $\Cap_p$-quasi continuous on $\Omega$.
\end{example}


\subsection{Capacitary Extremals}

In this subsection we assume that $X$ is a tms of type $\Theta$, $V$ is of class
$\Upsilon(X)$ and $\psi$ is $V$-admissible. Here we will prove existence and uniqueness
of capacitary extremals.

\begin{theorem}\label{thm:ypo}
  Let $A\subset X$ and $\su\in V$ be non-negative. Then $\su\in\overline\cY_V(A)$
  if and only if $\tilde\su\geq 1$ $\Cap_{\psi,V}$-q.e. on $A$.
\end{theorem}

\begin{proof}
  When $\tilde\su\geq 1$ $\Cap_{\psi,V}$-q.e.\ on $A$, then $\su\in\overline\cY_V(A)$ by
  Lemma \ref{lem:approx}. For the converse implication let $\su\in\overline\cY_V(A)$.
  By Theorem \ref{thm:convsubse} there exists a sequence
  $(\su_n)_n\subset\cY_V(A)$ such that $\tilde\su_n\to\tilde\su$
  $\Cap_{\psi,V}$-q.e. on $X$. For every $n\in\IN$ there exists
  an open set $O_n$ in $X$ containing $A$ such that $\su_n\geq 1$ $\mu$-a.e. on $O_n$.
  Hence $\tilde\su_n\geq 1$ $\Cap_{\psi,V}$-q.e.  on $A$ by Theorem \ref{thm:esti}. 
  This shows that $\tilde\su\geq 1$ $\Cap_{\psi,V}$-q.e. on $A$.
\end{proof}

\begin{theorem}\label{thm:quasi-cap}
  For $A\subset X$ the $\Upsilon$-capacity $\Cap_{\psi,V}$ of $A$ is given by
  \begin{eqnarray}
     \Cap_{\psi,V}(A) 
     &=& \inf\klg{\psi(\rho(\su)):\su\in\sV,\su\geq 1\text{ $\Cap_{\psi,V}$-q.e. on }A} \\
     \label{eq:inf}
     &=& \inf\klg{\psi(\rho(\su)):\su\in V,\tilde\su\geq 1\text{ $\Cap_{\psi,V}$-q.e. on }A}.
  \end{eqnarray}
\end{theorem}

\begin{proof}
  Let $I$ denote the infimum on the right hand side of \eqref{eq:inf} and $\su\in\cY_V(A)$.
  Then by Theorem \ref{thm:ypo} we get that $\tilde\su^+\geq 1$ $\Cap_{\psi,V}$-q.e. on $A$.
  Hence $I\leq \psi(\rho(\su^+))\leq\psi(\rho(\su))$. Taking the infimum over all
  $\su\in\cY_V(A)$ we get that $I\leq \Cap_{\psi,V}(A)$.
  On the other hand, let $\su\in V$ be such that $\tilde\su\geq 1$
  $\Cap_{\psi,V}$-q.e. on $A$. By Lemma \ref{lem:approx} there exist
  $\su_n\in\cY_V(A)$ such that $\su_n\to\su^+$ in $V$. Hence 
  $\Cap_{\psi,V}(A)\leq\psi(\rho(\su_n))\rightarrow\psi(\rho(\su^+))\leq\psi(\rho(\su))$.
  Taking the infimum over all such $\su$ gives that $\Cap_{\psi,V}(A)\leq I$
  and hence we have equality.
\end{proof}

\begin{definition}\label{def:capextremal}
  A function $\su\in\sV$ is called a/the {\em $\Cap_{\psi,V}$-extremal for
  $A\subset X$} if $\su\geq 1$ $\Cap_{\psi,V}$-q.e on $A$ and
  $\psi(\rho((\su))=\Cap_{\psi,V}(A)$.
\end{definition}

\begin{theorem}\label{thm:extremal}
  For every $A\subset X$ with $\Cap_{\psi,V}(A)<\infty$ there exists a
  $\Cap_{\psi,V}$-extremal $\se_A\in\sV$ with $0\leq\se_A\leq 1$ $\Cap_{\psi,V}$-q.e. on $X$
  and $\se_A=1$ $\Cap_{\psi,V}$-q.e. on $A$. If in addition
  $\rho$ is strictly convex, then $\se_A$ is unique.
\end{theorem}

\begin{proof}
  Let $\cY_V^+(A):=\klg{\su\in\cY_V(A):\su=\su^+}$.
  Since $\Cap_{\psi,V}(A)<\infty$ we have that $\cX:=\overline\cY^+_V(A)$ is a non-empty
  closed and convex subset of $V$. Let $(\su_n)_n\subset\cY_V^+(A)$ be such that
  $\psi(\rho(\su_n))\to\Cap_{\psi,V}(A)$. Then the sequence $(\su_n)_n$ is bounded
  in the reflexive Banach space $V$ and hence, by possibly passing to a subsequence,
  weakly convergent to a function $\su\in\cX$. Using the weak lower semi-continuity
  of $\rho$ we get that $\psi(\rho(\su))\leq \liminf_n \psi(\rho(\su_n))=\Cap_{\psi,V}(A)$.
  Since $\psi(\rho(\sv))\geq\Cap_{\psi,V}(A)$ for all $\sv\in\cY_V^+(A)$
  and hence for all $\sv\in\cX$ we get $\psi(\rho(\su))=\Cap_{\psi,V}(A)$.
  From Theorem \ref{thm:ypo} we get that $\se_A:=\tilde\su\geq 1$
  $\Cap_{\psi,V}$-q.e. on $A$. By possibly replacing $\se_A$ with $(se_A\wedge 1)^+$
  the existence part is proved. Uniqueness follows from the strict convexity of
  $\rho$ and from the fact that $\psi$ is strictly increasing.
\end{proof}

\begin{remark}\label{rem:inypo}
  The $\Cap_{p,\Omega}$-extremal for $A\subset\overline\Omega$ is the projection of
  $0$ onto $\overline\cY_{p,\Omega}(A)$.
\end{remark}

\section{Vanishing 'boundary' values}

In this section we give an application of the $\Upsilon$-capacity, namely
to decide if a given function $\su$ lies in $W^{1,\Phi}_0(\Omega)$ or not.
Here $\Phi\in\dt\cap\nt$ is an $\cN$-function and $W^{1,\Phi}_0(\Omega)$ is the
closure of $W^{1,\Phi}(\Omega)\cap C_c(\Omega)$ in $W^{1,\Phi}(\Omega)$.
We will assume in this section that $X$ is a tms of type $\Theta$ which
satisfies the {\em second axiom of countability}, $V$ is of class $\Upsilon(X)$
and $\psi$ is $V$-admissible and
\begin{itemize}
  \item[(V7)] $\su_n,\sv\in V$, $\su_n\to\su$ in $V$ implies that
              $\su_n\wedge\sv\to \su\wedge\sv$ in $V$.
\end{itemize}

\begin{definition}
  For $\su\in V$ we let $\supp(\su)$ be the intersection of all closed
  sets $A\subset X$ such that $\su=0$ $\mu$-a.e.\ on $X\setminus A$. Note that
  $O:=X\setminus\supp(\su)$ is the largest open set in $X$ such that
  $\su=0$ $\mu$-a.e. on $O$.
  For a set $X_0\subset X$ we let $V_c(X_0)$ be the space consisting
  of all $\su\in V$ such that $\supp(\su)\subset X_0$ is compact and we let
  $V_0(X_0)$ be the closure of $V_c(X_0)$ in $V$.
\end{definition}

\begin{theorem}
  Let $X_0\subset X$ be non-empty and assume that $V_\infty:=V\cap L^\infty(X,\mu)$ with
  respect to the norm $\norm{\su}_{V_\infty}:=\norm{\su}_V+\norm{\su}_{L^\infty}$ is a
  Banach algebra and $V$ has the cutoff-property. Then
  \begin{equation}\label{eq:v0}
     V_0(X_0) = \klgk{\su\in V:\tilde\su=0\;\Cap_{\psi,V}\text{-q.e.\ on } X\setminus X_0}.
  \end{equation}
\end{theorem}

\begin{proof}
  Let $D_0$ denote the right hand side of \eqref{eq:v0}. First we show that
  $V_0(X_0)\subset D_0$. Let $\su\in V_0(X_0)$. Then there exists a sequence
  of functions $\su_n\in V_c$ such that $\su_n\to \su$ in $V$.
  By possibly passing to a subsequence (Theorem \ref{thm:convsubse}) we get
  that $(\tilde\su_n)_n$ converges $\Cap_{\psi,V}$-quasi everywhere to $\tilde\su$ and
  hence (Theorem \ref{thm:esti}) $\tilde\su=0$ $\Cap_{\psi,V}$-quasi everywhere on
  $X\setminus X_0$, that is, $\su\in D_0$.

  To show that $D_0\subset V_0(X_0)$ we first consider a non-negative function
  $\su\in D_0\cap L^\infty(X)$. Then there exists a sequence $(u_n)_n$ in
  $V\cap C_c(X)$ which converges to $\su$ in $V$. Since
  $(u_n\vee 0)\wedge\norm{\su}_\infty$ converges also to $\su$ in $V$
  we may assume that $0\leq u_n\leq\norm{\su}_\infty$. Let $u\in\tilde\su$ be fixed. 
  By possibly passing to a subsequence (Theorem \ref{thm:convsubse})
  we get that for each $m\in\IN$ there exists an open set $G_m$ in $X$ such that
  $\Cap_{\psi,V}(G_m)\leq 1/m$ and $u_n\to u$ uniformly on $X\setminus G_m$.
  Hence there exists $n_0=n_0(m)$ such that $\betragk{u_{n_0}-u}\leq 1/(2m)$ everywhere on
  $X\setminus G_m$ and $\normk{u_{n_0}-\su}_V\leq 1/m$.
  Let $U_m$ be an open set in $X$ such that $\Cap_{\psi,V}(U_m)\leq 1/m$ and $u=0$
  everywhere on $X\setminus(X_0\cup U_m)$. Consequently, $\betragk{u_{n_0}}\leq 1/(2m)$ everywhere
  on $X\setminus (X_0\cup O_m)$ where $O_m:=G_m\cup U_m$.
  Let $\se_m\in\sV$ be a $\Cap_{\psi,V}$-extremal for $O_m$ (Theorem \ref{thm:extremal})
  and fix $e_m\in\se_m$. By possibly changing $e_m$ on a $\Cap_{\psi,V}$-polar set we
  may assume that $e_m\equiv 1$ everywhere on $O_m$ and $0\leq e_m\leq 1$ everywhere
  on $X$. For $w_m:=(u_{n_0}-1/m)^+$ we have that $v_m:=w_m(1-e_m)\in V_c$. In fact, 
  \[ \supp(v_m) 
      \subset \supp(w_m)\cap \supp(1-e_m)
      \subset \klg{x\in X:u_{n_0}(x)\geq 1/m}\cap O_m^c \subset (X_0\cup O_m)\cap O_m^c \subset X_0.
  \]
  Since $v_m=w_m-w_me_m$ is a bounded sequence in the reflexive Banach space
  $(V_0(X_0),\norm{\cdot}_V)$, by possibly passing to a subsequence, we may assume that
  $v_m\rightharpoonup \sv\in V_0(X_0)$. By Mazur's lemma
  [using that by possibly passing to a subsequence $v_m\to\su$ $\mu$-a.e.\ and property \prss]
  we get that $\sv=\su$ $\mu$-a.e.\ on $X$ and hence $\su\in V_0(X_0)$.
  If $\su\in D_0\cap L^\infty(X,\mu)$ is arbitrary, then we get by what we proved already
  that $\su^+$ and $\su^-$ belong to $V_0(X_0)$ and hence $\su\in V_0(X_0)$.
  Finally, if $\su\in D_0$, then there exist $u_n\in V\cap C_c(X)\subset L^\infty(X,\mu)$
  such that $u_n\to \su$ in $V$. Let $\sw_n:=(u_n\wedge\su^+)\vee (-\su^-)\in D_0\cap L^\infty(X,\mu)$.
  Then $\sw_n\in V_0(X_0)$ and hence $\su=\lim_n \sw_n\in V_0(X_0)$ [property (V7)].
\end{proof}

\begin{definition}
  Let $\Omega$ be an open and non-empty set in $\IR^N$, $p\in (1,\infty)$ and
  $\Phi\in\dt\cap\nt$ be an $\cN$-function. Then we define
  $W^{1,p}_0(\Omega)$ and $W^{1,\Phi}(\Omega)$ as the closure of $\cD(\Omega)$ in
  $W^{1,p}(\Omega)$ and $W^{1,\Phi}(\Omega)$, respectively.
\end{definition}

\begin{corollary}
  Let $\Omega$ be an open and non-empty set in $\IR^N$, $p\in (1,\infty)$ and
  $\Phi\in\dt\cap\nt$ be an $\cN$-function. Then
  $W^{1,p}_0(\Omega) = \klg{\su\in\tW^{1,p}(\Omega):\tilde\su=0\;
               \text{ $\Cap_{p,\Omega}$-q.e. on $\partial\Omega$} }$ and
    $W^{1,\Phi}_0(\Omega) = 
          \klg{\su\in\tW^{1,\Phi}(\Omega):\tilde\su=0\;
               \text{ $\Cap_{\Phi,\Omega}$-q.e. on $\partial\Omega$} }$.
\end{corollary}

To finish this section and the article we mention two further characterizations of $W^{1,p}_0(\Omega)$.
The original proof of Theorem \ref{thm:havin} is due to Havin \cite{havin:68:aim} and
Bagby \cite{bagby:72:qut}, an alternative proof is given by Hedberg \cite{hedberg:72:nlp}.
An other characterization, Theorem \ref{thm:ziemer}, was recently proved by David Swanson and
William P. Ziemer \cite[Theorem 2.2]{swanson:99:sfi}. The main difference to
Theorem \ref{thm:havin} is that the function $\su$ was not assumed to belong to the space $W^{1,p}(\IR^N)$.

\begin{theorem}\label{thm:havin}
  Let $1<p<\infty$, $\Omega\subset\IR^N$ an open set and let $\su\in W^{1,p}(\IR^N)$.
  Then $\su\in W^{1,p}_0(\Omega)$ if and only if
  $\lim_{r\to 0} r^{-N}\int_{B(x,r)} \betrag{u(y)}\;dy = 0$
  for $Cap_{p}$-q.e. $x\in\IR^N\setminus\Omega$. 
\end{theorem}

\begin{theorem}\label{thm:ziemer}.
  Let $p\in(1,\infty)$ and $\su\in W^{1,p}(\Omega)$. If
  $\lim_{r\to 0} r^{-N} \int_{\cB(x,r)\cap\Omega} \betrag{\su(y)}\;dy = 0$
  for $\Cap_p$-quasi every $x\in\partial\Omega$, then $\su\in W^{1,p}_0(\Omega)$.
\end{theorem}

\bibliography{biblio2}

\providecommand{\mathbb}[1]{\mathbf{#1}}\providecommand{\cprime}{$'$}
\providecommand{\bysame}{\leavevmode\hbox to3em{\hrulefill}\thinspace}
\providecommand{\MR}{\relax\ifhmode\unskip\space\fi MR }
\providecommand{\MRhref}[2]{%
  \href{http://www.ams.org/mathscinet-getitem?mr=#1}{#2}
}
\providecommand{\href}[2]{#2}
\begin{thebibliography}{10}

\bibitem{adams:96:fsp}
David~R. Adams and Lars~Inge Hedberg, \emph{Function spaces and potential
  theory}, Grundlehren der mathematischen Wissenschaften, vol. 314,
  Springer-Verlag, Berlin, 1996. \MR{97j:46024}

\bibitem{arendt:03:lrb}
Wolfgang Arendt and Mahamadi Warma, \emph{The {L}aplacian with {R}obin boundary
  conditions on arbitrary domains}, Potential Anal. \textbf{19} (2003), no.~4,
  341--363. \MR{1 988 110}

\bibitem{bagby:72:qut}
Thomas Bagby, \emph{{Quasi topologies and rational approximation.}}, J. Funct.
  Anal. \textbf{10} (1972), 259--268.

\bibitem{bouleau:91:df}
Nicolas Bouleau and Francis Hirsch, \emph{Dirichlet forms and analysis on
  {W}iener space}, de Gruyter Studies in Mathematics, vol.~14, Walter de
  Gruyter \& Co., Berlin, 1991. \MR{MR1133391 (93e:60107)}

\bibitem{brezis:83:af}
Ha{\"{\i}}m Brezis, \emph{Analyse fonctionnelle}, Collection Math\'ematiques
  Appliqu\'ees pour la Ma\^\i trise. [Collection of Applied Mathematics for the
  Master's Degree], Masson, Paris, 1983, Th{\'e}orie et applications. [Theory
  and applications]. \MR{MR697382 (85a:46001)}

\bibitem{choquet:54:toc}
Gustave Choquet, \emph{{Theory of capacities.}}, Ann. Inst. Fourier (Grenoble)
  \textbf{5} (1954), 131--295.

\bibitem{evans:92:mtf}
Lawrence~C. Evans and Ronald~F. Gariepy, \emph{Measure theory and fine
  properties of functions}, Studies in Advanced Mathematics, CRC Press, Boca
  Raton, FL, 1992. \MR{MR1158660 (93f:28001)}

\bibitem{koskela:08:see}
Piotr Haj{\l}asz, Pekka Koskela, and Heli Tuominen, \emph{Sobolev embeddings,
  extensions and measure density condition}, J. Funct. Anal. \textbf{254}
  (2008), no.~5, 1217--1234. \MR{MR2386936}

\bibitem{havin:68:aim}
Victor~P. Havin, \emph{Approximation in the mean by analytic functions}, Dokl.
  Akad. Nauk SSSR 9 (1968), 245--248.

\bibitem{hedberg:72:nlp}
Lars~Inge Hedberg, \emph{{Non-linear potentials and approximation in the mean
  by analytic functions.}}, Math. Z. \textbf{129} (1972), 299--319.

\bibitem{heinonen:93:npt}
Juha Heinonen, Tero Kilpel{\"a}inen, and Olli Martio, \emph{Nonlinear potential
  theory of degenerate elliptic equations}, Oxford Mathematical Monographs,
  Clarendon Press, New York, 1993. \MR{94e:31003}

\bibitem{hudzik:00:npmo}
Henryk Hudzik and Pawe{\l} Kolwicz, \emph{A note on {$P$}-convexity of {O}rlicz
  and {M}usielak-{O}rlicz spaces of {B}ochner type}, Function spaces
  ({P}ozna\'n, 1998), Lecture Notes in Pure and Appl. Math., vol. 213, Dekker,
  New York, 2000, pp.~181--191. \MR{MR1772124 (2001i:46060)}

\bibitem{ma:09:fbfos}
J.~Mal\'y, D.~Swanson, and W.~P. Ziemer, \emph{Fine behavior of functions whose
  gradients are in an orlicz space}, Studia Math. \textbf{190(1)} (2009),
  33--71.

\bibitem{ziemer:97:fr}
Jan Mal{\'y} and William~P. Ziemer, \emph{Fine regularity of solutions of
  elliptic partial differential equations}, Mathematical Surveys and
  Monographs, vol.~51, American Mathematical Society, Providence, RI, 1997.
  \MR{MR1461542 (98h:35080)}

\bibitem{mazya:85:ssp}
Vladimir~G. Maz{'}ya, \emph{Sobolev spaces}, Springer Series in Soviet
  Mathematics, Springer-Verlag, Berlin, 1985, Translated from the Russian by T.
  O. Shaposhnikova. \MR{87g:46056}

\bibitem{mu:83:osms}
Julian Musielak, \emph{Orlicz spaces and modular spaces}, Lecture Notes in
  Mathematics, vol. 1034, Springer-Verlag, Berlin, 1983. \MR{MR724434
  (85m:46028)}

\bibitem{rr:91:toos}
M.~M. Rao and Z.~D. Ren, \emph{Theory of {O}rlicz spaces}, Monographs and
  Textbooks in Pure and Applied Mathematics, vol. 146, Marcel Dekker Inc., New
  York, 1991. \MR{MR1113700 (92e:46059)}

\bibitem{rr:02:aoos}
\bysame, \emph{Applications of {O}rlicz spaces}, Monographs and Textbooks in
  Pure and Applied Mathematics, vol. 250, Marcel Dekker Inc., New York, 2002.
  \MR{MR1890178 (2003e:46041)}

\bibitem{shvartsman:07:eos}
P.~Shvartsman, \emph{On extensions of {S}obolev functions defined on regular
  subsets of metric measure spaces}, J. Approx. Theory \textbf{144} (2007),
  no.~2, 139--161. \MR{MR2293385 (2007k:46057)}

\bibitem{swanson:99:sfi}
David Swanson and William~P. Ziemer, \emph{Sobolev functions whose inner trace
  at the boundary is zero}, Ark. Mat. \textbf{37} (1999), no.~2, 373--380.
  \MR{2000g:46048}

\end{thebibliography}
\bibliographystyle{amsplain}

\end{document}